\documentclass{amsart}
\input epsf
\newtheorem{theorem}[equation]{Theorem}

\newtheorem{lemma}[equation]{Lemma}

\theoremstyle{remark}
\newtheorem{remark}[equation]{Remark}
\theoremstyle{definition}

\numberwithin{equation}{section}
\renewcommand{\qed}{\hspace*{\fill} \setlength{\unitlength}{1mm}
\begin{picture}(2.5,2.5)
      \put(0,0){\framebox(2.5,2.5){}}
  \end{picture}
\setlength{\unitlength}{1pt}}

\newcommand{\reals}{{\bf R}}

\newcommand{\cN}{{\mathcal{N}}}

\def\bC{{\mathbf{C}}}

\begin{document}
\title{On nodal sets and nodal domains on $S^2$ and $\reals^2$}
\author{Alexandre Eremenko, Dmitry Jakobson and Nikolai Nadirashvili}
\address{Mathematics Department, Purdue University,
150 N University Street W. Lafayette, IN 47907-2067 USA; e-mail
\texttt{eremenko@math.purdue.edu}
\newline Department of Mathematics and Statistics, McGill University, 805
Sherbrooke Str. West, Montreal, QC H3A 2K6, Canada; e-mail
\texttt{jakobson@math.mcgill.ca} \newline Laboratoire d'Analyse,
Topologie, Probabilit\'{e}s UMR 6632, Centre de Math\'{e}matiques et
Informatique, Universit\'{e} de Provence, 39 rue F. Joliot-Curie,
13453 Marseille Cedex 13, France; e-mail
\texttt{nicolas@cmi.univ-mrs.fr} }

\thanks{The first author was supported by NSF grants
DMS-0555279 and DMS-244547.  The second author was supported by
NSERC, FQRNT and Dawson fellowship.}

\subjclass[2000]{Primary: 58J50 Secondary: 11J70, 35P20, 81Q50}

\date{\today}
\begin{abstract}
We discuss possible topological configurations of nodal sets, in
particular the number of their components,
for
spherical harmonics on $S^2$.  We also construct a solution of the
equation $\Delta u=u$ in $\reals^2$
that has only two nodal
domains. This equation arises in the study of high energy eigenfunctions.
\end{abstract}

\maketitle

\section{Topological structure of nodal domains}
Homogeneous harmonic polynomials of degree $n$ in three variables,
when restricted to the unit sphere, become eigenfunctions of the
Laplace--Beltrami operator on the sphere, with eigenvalue
$\lambda=n(n+1)$ of multiplicity $2n+1$, $n=0,1,\ldots.$ We
permit ourselves the liberty of calling them {\em eigenfunctions of
degree $n$}.

The zero set of an eigenfunction is called the {\em nodal set}. The nodal
set is {\em non-singular} if it is a union of disjoint closed
analytic curves, and the nodal set of a generic eigenfunction is
non-singular \cite{generic}.

The known results about topology of nodal sets on the sphere are the
following \cite{Survey}. The nodal set is non-empty for $n\geq 1$. A
general theorem of Courant \cite{CH} implies that the nodal set of
an eigenfunction of degree $n$ consists of at most $n^2$ components.
H. Lewy \cite{Lewy} proved that for even $n\geq 2$ the number of
components is at least $2$, and constructed eigenfunctions of any
degree $n$ whose nodal sets have one component for odd $n$ and two
components for even $n\geq 2$.

We call subsets $X_1$ and $X_2$ of homeomorphic topological spaces
$U_1$ and $U_2$ {\em equivalent} if there is a homeomorphism
$h:U_1\to U_2$ such that $h(X_1)=X_2$. The ambient spaces $X_1$ and $X_2$
are usually clear from context. We notice that each
homogeneous harmonic polynomial is either even or odd, so its zero
set is invariant under the antipodal map of the sphere.
\vspace{.1in}

\begin{theorem}\label{thm1}
Let $0<m\leq n$, and let $n-m$ be even. For every set of $m$
disjoint closed curves on the sphere, whose union $E$ is invariant
with respect to the antipodal map, there exists an eigenfunction of
degree $n$ whose zero set is equivalent to $E$.
\end{theorem}

The invariance with respect to the antipodal map $a$ is a strong
restriction on the possible shape of the set (see, for example,
\cite{Viro} for the proofs of the following facts). For each
component $C$ of the nodal set, either $a(C)=C$, or $a(C)=C'$, where
$C'$ is another component, different from $C$. In the former
case $C$ is called an {\em odd component}, and in the latter case an
{\em even component} or an {\em oval}. Odd eigenfunctions have
exactly one odd component of the nodal set, and even eigenfunctions
have only even components in their nodal sets.

Thus the total number of components is of the same
parity as $n$. This generalizes Lewy's result that
there are at least two components for even $n\geq 2$.

All these facts hold for {\em any} homogeneous polynomials
in three variables restricted to the sphere with non-singular
zero sets.

Topological classification of non-singular level sets of arbitrary
(not necessarily harmonic) homogeneous polynomials was obtained in
\cite{santos}; the author also allows some simple singularities
(double points).

\begin{remark}\label{remark:upbound}
Theorem \ref{thm1} gives at most $n$ components of the nodal set for
a spherical harmonic of degree $n$, while the best upper estimates
known to the authors are $n^2+O(n)$.
These estimates can be derived from the upper
estimates of the number of nodal
domains,
\cite{Karp02} (see also
\cite[Problem 1]{unsolved}, and \cite{Karp89,Karp94}). Thus we obtain
from the Corollary in
\cite{Karp02} that the number of components of
the nodal set of a spherical harmonic of
degree $n$ does not exceed
\begin{itemize}
\item[i)]
$n^2-2n+2$ for even $n$, and
\item[ii)]
$(n-1)^2+3$ for odd $n$.
\end{itemize}
For $n\leq 6$, precise estimates of the number of nodal domains
were obtained
in \cite{Leydold}.
For large $n$, Pleijel \cite{Pleijel} obtained an upper bound
$(4+o(1))n^2/j_0^2\approx 0.69n^2$,
where $j_0$ denotes the smallest zero of the
$0$-th Bessel function.
   It seems difficult to determine the largest
possible number of components of the nodal set of a spherical
harmonic of given degree. We
address this question in section \ref{sec:ovals}.
\end{remark}

In the proof of Theorem \ref{thm1} we use a related result on
topological classification of zero sets of harmonic polynomials in
two variables (Theorem \ref{thm2} below). To state Theorem
\ref{thm2} we need to fix some terminology.

%%This result is essentially well-known though we could not find its
%%statement in the literature.

A {\em tree} is a finite connected contractible $1$-complex with at
least one edge. Vertices of degree $1$ are called {\em leaves}. A
{\em forest} is a disjoint union of trees. An {\em embedded forest}
is a subset of the plane which is the image of a proper embedding of
a forest minus leaves to the plane.

Let $u$ be a harmonic polynomial of two variables of degree $n$.
Then it is easy to see that the level set $\{ z:u(z)=0\}$ is an
embedded forest (indeed, existence of a cycle would contradict the
maximum principle).  All vertices of a forest, except the leaves,
 have even
degrees (since the function changes sign an even number of times as
we go around the vertex); there are exactly $2n$ leaves.

\vspace{.1in}

\begin{theorem}\label{thm2}
Let $F$ be an embedded forest with $2n$ leaves and such that all its
vertices in the plane are of even degrees. Then there exists a
harmonic polynomial $u$ of degree $n$ whose zero set is equivalent
to $F$.
\end{theorem}

To prove Theorem \ref{thm1}, we need only the ``generic case'' of
Theorem \ref{thm2}, when $F$ is a union of simple curves, that is
each tree of the forest has only one edge. Such forest is convenient
to visualize as a {\em chord diagram} in the unit disc. It is
obtained by an embedding of a forest to the closed unit disc such
that the leaves are mapped to the boundary circle and the rest of
the forest into the open disc. It is convenient to place the leaves
at the roots of unity of degree $2n$ on the unit circle which we
will always do. \vspace{.1in}

\begin{remark} The generic case of Theorem \ref{thm2} (needed in the
proof of Theorem \ref{thm1}) can be derived from Belyi's theorem, see e.g.
\cite[Theorem 2.2.9]{LZ}. Using the full strength of Theorem
\ref{thm2}, and some extra work, one can prove an extended version
of Theorem \ref{thm1}, classifying all nodal sets, not only generic
ones. Namely, one can show that any gluing of a forest with its
image under the antipodal map (gluing by identifying the leafs in the
natural cyclic order) is
equivalent to a nodal set, and that all nodal sets arise this way.
\end{remark}

The connection between the level sets of harmonic polynomials in
the plane and the nodal sets of eigenfunctions is established by the
following Lemma~\ref{levelset}.

We recall that every eigenfunction
of degree $n$ can be represented in the upper hemisphere as
$$f(z)=\Re \sum_{k=0}^n L_k(r)a_kz^k,$$
where $z=x+iy$, is the orthogonal projection of a point of the upper
hemisphere onto the equatorial plane which we identify with the
complex plane, $r=|z|$ and $L_k(r)=F_n^{(k)}(\sqrt{1-r^2}),$ where
$F_n^k$ are defined by
\begin{equation}\label{assoc:legendre}
c_{n,k}(1-x^2)^{k/2}F_n^k(x)=\frac{d^k}{dx^k}P_n(x).
\end{equation}
Here $P_n$ is the Legendre polynomial of degree $n$ (see, for
example, \cite{TS}):
\begin{equation}\label{legendre}
P_n(x)=\frac{1}{2^n n!}\;\frac{d^n}{dx^n}\left([x^2-1]^n\right),
\end{equation}
and where $c_{n,k}$ is chosen so that
\begin{equation}\label{value:legendre}
L_k(0)=F_n^k(1)=1.
\end{equation}

Let
$$p(z)=a_0+\ldots+z^n$$
be a monic polynomial whose real part has non-singular zero set.
We define a family of eigenfunctions
$$
f_t(z)=\Re \sum_{k=0}^n L_k(r)t^{n-k} a_kz^k,
$$
where $t>0$ is a parameter.
\vspace{.1in}

\begin{lemma}\label{levelset}
If $t$ is small enough then the zero set
of $f_t$ in the open upper hemisphere is equivalent to the
zero set of $\Re p$ in the plane.
\end{lemma}

{\em Proof.} We identify the upper hemisphere with the unit disc in
the complex plane via vertical projection as before. Consider the
transformed function
$$t^{-n}f_t(t z)=
\Re \sum_{k=0}^n L_k(t r)a_kz^k.$$

When $t\to 0$ this converges to $p(z)$ in $C^1(K)$ on every
disc $K$ in the plane. Here we used the property
\eqref{value:legendre} of Legendre polynomials. Thus for every
$r_0>0$ there exists $t_0>0$ such that for
$t<t_0$,
the zero set of $t^{-n}f_t(t
z)$ is equivalent to the zero set of $\Re p$ in the disc $\{ z:|z|<r_0\}$.
We fix
$r_0$ in such a way that the portion of the zero set of $\Re p$ in
$|z|<r_0$ is equivalent to its zero set in the whole plane, and thus
the portion of the zero set of $\Re p$ in in $r_0\leq |z|<\infty$ is
equivalent to the zero set of $\Re z^n$.

$C^1$-convergence implies that the zero set of $f_t$ with
$t<t_0$ in the disc $|z|<r_0$
is equivalent to
the level set of $p$ in the plane. Now we consider the set
$r_0\leq |z|\leq 1$. It is clear that $f_t(z)\to z^n$
in $C^1$ on this set, as $t\to 0$. So for $t$ small
enough the zero set of $f_t$ in the unit disc is equivalent to
the zero set of $p$ in the plane. This proves the lemma. \qed

\vspace{.1in}

{\em Proof of Theorem \ref{thm1}}. In view of Lemma \ref{levelset},
it remains to show that for every integer $m\in [1,n]$ of the same
parity as $n$, every disjoint union $E$ of $m$ closed curves on the
sphere that is invariant with respect to the antipodal map can be
obtained by gluing together a chord diagram and its centrally
symmetric chord diagram.
\vspace{.1in}

{\em Case of odd $m$ and $n$, $m=2k+1$}. As we stated after Theorem
\ref{thm1}, there is exactly one odd component $C_0$ (mapped into
itself by the antipodal map).  Let $C_1,C_1',\ldots,C_k,C_k'$ be the
other components, such that $C_j$ is mapped into $C_j'$ by the antipodal
map $T$.

First, remove the odd component $C_0$ from the sphere. Its
complement consists of two disks $D$ and $D'$, and we
can choose our notation so that $C_j\subset D$
and $C_j^\prime\subset D^\prime$ for $1\leq j\leq k$.  $T$ induces a
homeomorphism between the two disks, while their common boundary
$C_0$ satisfies $T(C_0)=C_0$.  Choose two antipodal points $a,b\in
C_0$.

Choose a simple curve $\gamma$ lying in $D$, so that it
\begin{itemize}
\item[i)] connects $a$ and $b$;
\item[ii)] intersects each $C_j$ transversally an even number of
times;
\item[iii)] the total number of intersections of $\gamma$ with
$\cup_{j=1}^kC_j$ is equal to $n-1$.
\end{itemize}
The case $n=m=2k+1$, when $\gamma$ intersects
each $C_j$ twice, is illustrated in Figure 1 (in the right part
of this figure,
$\gamma$ is the upper half of the dotted circle).
If $m<n$, we can construct such a
$\gamma$ as follows: first choose a curve $\widetilde{\gamma}$
that intersects each $C_j$ exactly twice, for a total of $2k<n-1$
intersections.  Then adjust $\widetilde{\gamma}$ so that it
intersects $C_1$ at two more points, and repeat the
procedure $(n-1)/2-k$ times until the number of intersections is
equal to $n-1$.

\begin{center}
\epsfxsize=4.0in \centerline{\epsffile{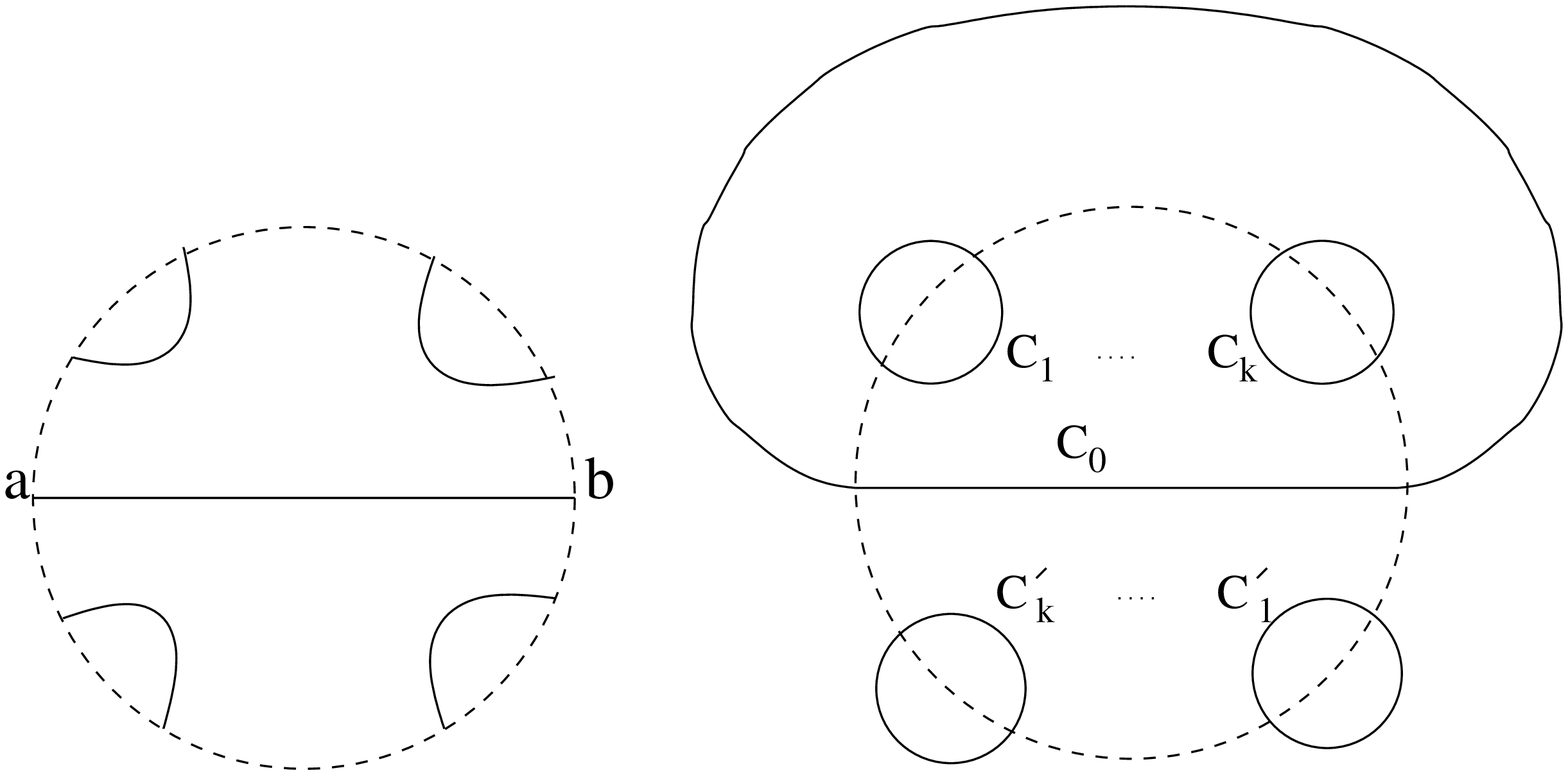}} Figure 1. A
chord diagram producing ovals.
\end{center}

Let $T(\gamma)=\gamma'$.  Then $\gamma\cup\gamma'\cup\{a\}\cup\{b\}$
is a simple closed curve on $S^2$ (shows as a dotted circle in Fig. 1).
Its complement consists of two
regions $G$ and $G'$ satisfying $T(G)=G'$.

Let us represent $G$ as a unit disk $\{|z|<1\}$, so that $a=1,b=-1$,
and $T(z)=-1/\bar{z}$.  The parts of curves $C_j$ and $C_j'$ lying
inside $G$, together with the arc of $C_0$ connecting $a$ and $b$
inside $G$, form a chord diagram that has $2n-2+2=2n$ vertices.

The parts of curves $C_j$ and $C_j'$ lying inside $G'$, together
with the complementary arc of $C_0$ connecting $a$ and $b$ inside
$G'$, form {\em another} chord diagram, which is symmetric to the
first diagram under the mapping $T(z)=-1/\bar{z}$.

An application of Lemma \ref{levelset} finishes the proof.
\vspace{.1in}

{\em Case of even $n$ and $m=2k$.} In this case all components
of the nodal set are even, and there is exactly one complementary
region to our set $E$ which is mapped by the antipodal map onto
itself. One can find a closed curve $\alpha$ in this region so that
$\alpha$ is invariant under the antipodal map. Clearly this $\alpha$
is disjoint from $E$. Adding it to $E$ we can perform the
construction as above, and then remove one edge corresponding to
$\alpha$ from the chord diagram. This completes the proof of Theorem
\ref{thm1}.
\qed
\vspace{.1in}

{\em Proof of Theorem \ref{thm2}}. We use the method of \cite{EG}.
A component of the complement of the embedded forest is called
a {\em face}. All faces are simply connected unbounded regions
The degree of a face is defined as the number of connected
components of its boundary. It is evident that the sum of the degrees of
all faces is $2n$.

A {\em labeling} of a forest is a prescribing of positive numbers to
the edges, so that the following condition is satisfied: for each
face, the sum of these numbers over all edges on the boundary of
this face equals $2\pi$ times the degree of the face.

It is easy to show that for every forest there exists a labeling. In
the special case that there are no vertices (the only case needed in
this paper) the labeling is unique and the label of each edge is
$2\pi$.

To define a labeling in the general case,
we notice that every forest can be obtained from a forest without
vertices by applying the following operation finitely many times.
Choose an edge $E=[a,\infty)$ whose one extremity is $\infty$.
Choose a point $v$ in the interior of this edge.
Then replace the part $[v,\infty)$ of $E$ by an odd number of edges
$E_1,\ldots,E_{2k+1}$ whose extremities are $v$ and $\infty$.
Thus $E$ is replaced by $E_0=[a,v]\subset E$ and $E_1,\ldots,E_{2k+1}$.
We suppose that $E_0,E_1,\ldots,E_{2k+1}$ are enumerated in a natural
cyclic order around $v$.

Now suppose that the label of $E$ was $x$. In the new tree, we
label $E_0$ by $x/2$, and then the labels
of $E_1,\ldots,E_{2k+1}$ will be $x/2$ and $2\pi-x/2$, alternatively.
It is clear that this gives a labeling of the new tree.

Once a labeling is prescribed, we define the ``length'' of a path in
the forest as the sum of the labels of edges of this path.

Now we orient the edges observing the following rule:
if $v$ is a vertex and $e_1,\ldots,e_{2k}$
are all edges attached to it, listed in the natural
cyclic order
then of any pair $e_j$ and $e_{j+1}$
one edge is oriented towards $v$ and another away from $v$.

It is easy to see that there are exactly two ways to
orient the edges of a tree according to this rule;
orientation of one edge induces orientation of the rest.
We begin by choosing orientation of the edges of one tree of
the forest. This induces orientation on the boundaries
of each face having common boundary with this tree.
So we obtain orientation of all other trees which intersect
the boundary of these faces. Continuing this procedure
we orient the edges of the whole forest, and in particular
the boundaries of all faces.

Now we construct a ramified covering from the forest
to the real line $\reals$. We equip the real line with
the spherical metric $2|dx|/(1+|x|^2)$, so that the length
of the real line is $2\pi$. We orient
the real line from $-\infty$ to $\infty$. Every edge of the
forest will be mapped homeomorphically and respecting
the orientation (that is increasing with respect to
the orientations of the edge and of the real line)
onto an interval
of $\reals$ whose length is the label of the edge.
All leaves are mapped to $\infty$.

Such a continuous map $\phi$ is uniquely defined once the
orientations of the edges and their labels are fixed.
It is a ramified covering (ramified at the vertices).

Now we extend $\phi$ to the faces.
Every boundary component of a face is mapped on the real
line homeomorphically. Our choice of the orientation of each
tree guarantees that the whole boundary (in $\bC$)
of each face is mapped on the real line by
a covering map. Considering the boundary of the face in $\bC$ and the real
line as circles we obtain a covering of circles.

We extend it to a ramified covering of the discs,
for example, with at most one ramification point in each disc
(none, if the boundary map is a homeomorphism).

Thus we obtain a ramified covering of the sphere,
so that the preimage of the real line is our embedded
forest, and the full preimage of infinity is infinity.

We pull back the standard conformal structure via $\phi$ and obtain
some conformal structure in the domain of $\phi$. By the
Uniformization theorem, there exists a homeomorphism $\psi$ such
that $\phi\circ\psi$ is holomorphic. As it is also of degree $n$
(=half of the number of leaves of the tree), we obtain a complex
analytic polynomial of degree $n$. The preimage of the real line is
evidently equivalent to our given forest.
\qed

\section{Number of ovals}\label{sec:ovals}

In this section for every $n\gg 0$, we construct examples of
spherical harmonics with many ovals. All spherical harmonics in this
section are assumed to be nonsingular; their nodal lines don't
intersect. Recall that a standard basis of spherical harmonics of
degree $n$ is given (up to proportionality constant which is
unimportant for us) by functions
$$
Y_n^m(\theta,\phi)=\sin^m\theta F_n^m(\cos\theta)\sin(m\phi), \qquad
0\leq m\leq n,
$$
where $F_n^m$ was defined in \eqref{assoc:legendre}. The function
$(1-x^2)^{m/2}F_n^m(x)$ is proportional to the {\em associated
Legendre function}; its zeros (aside from $x=\pm 1$) coincide with
those of $F_n^m(x)$.

For a spherical harmonic $\Psi$, we denote by ${\rm ov}(\Psi)$
the number of disjoint connected components (``ovals'') of its nodal
set.

We prove the following
\begin{theorem}\label{thm:ovals}
For any $n\gg 0$, there exist $\epsilon_n>0$ and a rotation
$R_n\in{\rm SO}(3)$ such that for any $\epsilon<\epsilon_n$,

\begin{itemize}
\item[i)] For {\em odd} $n$,
$$
\lim_{n\to\infty}\frac{{\rm ov}(Y_n^{[n/2]}+\epsilon\cdot
Y_n^{n-1}\circ R_n)}{n^2/4}=1.
$$
\item[ii)] For {\em even} $n$,
$$
\lim_{n\to\infty}\frac{{\rm ov}(Y_n^{n/2}+\epsilon\cdot Y_n^n\circ
R_n)}{n^2/4}=1.
$$
\end{itemize}

\end{theorem}

\begin{remark}
Theorem \ref{thm:ovals} complements results in \cite{Lewy}, where
spherical harmonics with the {\em smallest} possible number of ovals
(and nodal domains) are constructed. Leydold conjectured in
\cite{Leydold} that spherical harmonics $Y_n^{[n/2]}$ (that we are
perturbing) have the {\em largest} possible number of nodal domains.
However, their nodal sets are connected.  {\em Upper} bounds on the
number of nodal domains for spherical harmonics were discussed in
Remark \ref{remark:upbound}; for {\em irreducible} curves Harnack's
bound gives better results, see \cite{Gudkov}. It seems interesting
to construct explicit examples of spherical harmonics whose nodal
sets consist of a {\em prescribed} number of ovals between one (or
two in the even degree case) and $n^2/4$.
\end{remark}

{\em Proof of Theorem \ref{thm:ovals}.}
For convenience, we give separate proofs for different $n ({\rm mod} 4)$.

{\em Case i): $n=4k+3$.} This is the easiest case.

The nodal sets of the harmonic $Y_n^{2k+1}$ consists of
\begin{itemize}
\item[i)] north and south poles $(\theta=0,\theta=\pi)$.
\item[ii)] $2k+2$ parallels $\theta=\theta_j,\theta=\pi-\theta_j$, where
$$
0<\theta_1<\theta_2<\ldots<\theta_{k+1}<\pi/2;
$$
here $F_n^{2k+1}(\cos\theta_j)=0$.
\item[iii)] $4k+2$ half-meridians $\phi=\pi j/(2k+1),0\leq j<4k+2$.
\end{itemize}
The nodal lines intersect at the poles; as well as at $(2k+2)(4k+2)$
points in the region $\{0<\theta<\pi\}$, where two nodal lines intersect
at each point.  We denote the intersection points in the upper hemisphere
$\{0<\theta<\pi/2\}$ by $P_j,1\leq j\leq (k+1)(4k+2)$; and the
intersection points in the lower hemisphere
$\{\pi/2<\theta<\pi\}$ by $Q_j,1\leq j\leq (k+1)(4k+2)$.

The nodal set of the harmonic $Y_n^{4k+2}$ consists of
\begin{itemize}
\item[i)] north and south poles $(\theta=0,\theta=\pi)$.
\item[ii)] the equator $\theta=\pi/2$;
\item[iii)] $8k+4$ half-meridians $\phi=\pi j/(4k+2),0\leq j<8k+4$.
\end{itemize}
The nodal domains of $Y_n^{4k+2}$ are $8k+4$ sectors in the upper
hemisphere: $\{0<\theta<\pi/2; 0<\phi<\pi/(4k+2)\}$,
$\{0<\theta<\pi/2; \pi/(4k+2)<\phi<2\pi/(4k+2)\}$, etc.
and symmetric $8k+4$ sectors in the lower hemisphere.

Choose $\tilde{R}_n$ to be the rotation with $NS$-axis, and with an
angle $\psi_n$ satisfying $0<-\psi_n<\pi/(4k+2)$.  Then it
is easy to see that
\begin{lemma}\label{lemma:3mod4}
There are exactly {\em two} nodal meridians of
$Y_n^{4k+2}\circ \tilde{R}_n$ lying
between any two nodal meridians of $Y_n^{2k+1}$.  Accordingly,
$Y_n^{4k+2}\circ \tilde{R}_n$ takes {\em negative} values at all the points
$P_j$; and {\em positive} values at all the points
$Q_j,\; 1\leq j\leq (k+1)(4k+2)$.
\end{lemma}

We can next adjust $\tilde{R}_n$ slightly to get a rotation $R_n$ so
that the conclusion of Lemma \ref{lemma:3mod4} still holds, and also
$Y_n^{4k+2}\circ R_n$ takes {\em positive} value at $S$ and {\em
negative} value at $N$.

It is then clear that if we choose $\epsilon_n$ small enough,
then the nodal set of $Y_n^{2k+1}+\epsilon_n\cdot Y_n^{4k+2}\circ R_n$
will consist of
\begin{itemize}
\item[i)] one ``long'' oval along the equator.
\item[ii)] $(k+1)(2k+1)$ ``small'' ovals bounding {\em positive} nodal
domains in the upper hemisphere.
\item[iii)] $(k+1)(2k+1)$ ``small'' ovals bounding {\em negative} nodal
domains in the lower hemisphere.
\end{itemize}

Altogether, the nodal set consists of $(k+1)(4k+2)+1\sim n^2/4$
disjoint ovals, finishing the proof for $n=4k+3$.

\begin{center}
\epsfxsize=4.0in \centerline{\epsffile{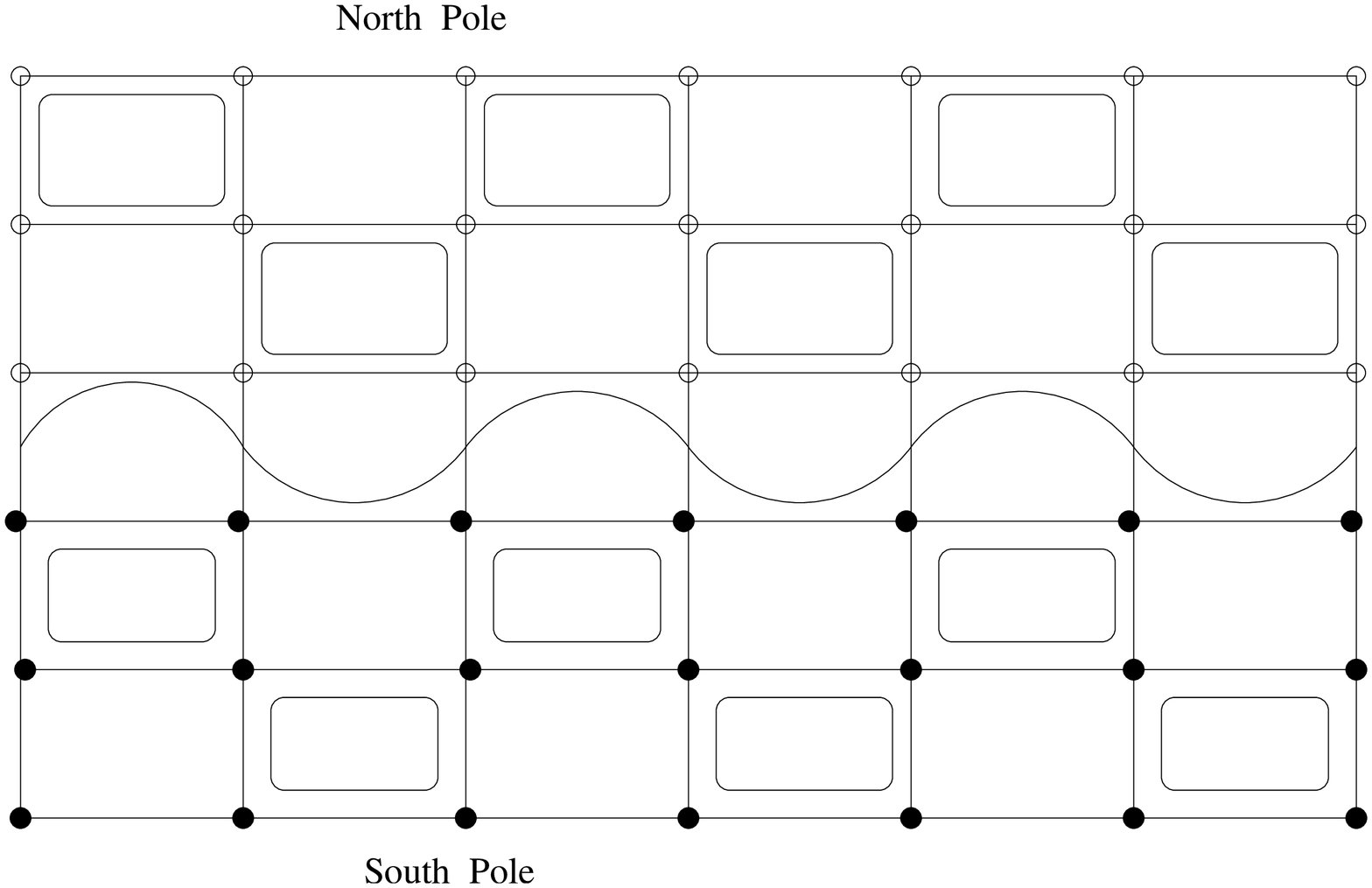}}

Figure 2. A perturbation in degree $4k+3=7$;\ $\bullet$\ denotes
positive sign of $Y_7^6$ at singular points of $Y_7^3$, while \
$\circ$\ denotes negative sign of $Y_7^6$ at singular points of
$Y_7^3$. Note that on $S^2$ all points on the top line (north pole)
and the bottom line (south pole) are identified.
\end{center}
\qed

%%%%%%%%%%%%%%%%

{\em Case ii): $n=4k+1$.}
The proof in this case proceeds along the same lines as in i),
except we lose a linear number of ovals due to the fact that it is
impossible to arrange for the perturbing harmonic to have constant
sign at so many adjacent nodal crossings as for $n=4k+3$.

Consider the spherical harmonic
$Y_{4k+1}^{2k}$.  Its nodal set consists of
\begin{itemize}
\item[i)] north and south poles $(\theta=0,\theta=\pi)$.
\item[ii)] $2k$ parallels $\theta=\theta_j,\theta=\pi-\theta_j$, where
$$
0<\theta_1<\theta_2<\ldots<\theta_k<\pi/2;
$$
here $F_n^{2k}(\cos\theta_j)=0$.
\item[iii)] $4k$ half-meridians $\phi=\pi j/(2k),0\leq j<4k$.
\end{itemize}
The nodal lines intersect at the poles; $4k$ points on the equator;
as well as at $(2k)(4k)$
points in the region $\{0<\theta<\pi\}$, where two nodal lines intersect
at each point.  We denote the intersection points in the upper hemisphere
$\{0<\theta<\pi/2\}$ by $P_j,1\leq j\leq k\cdot 4k$; and the
intersection points in the lower hemisphere
$\{\pi/2<\theta<\pi\}$ by $Q_j,1\leq j\leq k\cdot 4k$.

The nodal set of the harmonic $Y_n^{4k}$ consists of
\begin{itemize}
\item[i)] north and south poles $(\theta=0,\theta=\pi)$.
\item[ii)] the equator $\theta=\pi/2$;
\item[iii)] $8k$ half-meridians $\phi=\pi j/(4k),0\leq j<8k$.
\end{itemize}
The nodal domains of $Y_n^{4k}$ are $8k$ sectors in the upper
hemisphere: $\{0<\theta<\pi/2; 0<\phi<\pi/(4k)\}$,
$\{0<\theta<\pi/2; \pi/(4k)<\phi<2\pi/(4k)\}$, etc.
and symmetric $8k$ sectors in the lower hemisphere.

Choose $\tilde{R}_n$ to be the rotation with $NS$-axis, and with an
angle $\psi_n$ satisfying $0<-\psi_n<\pi/(4k)$.  Then it
is easy to see that
\begin{lemma}\label{lemma:1mod4}
There are exactly {\em two} nodal meridians of
$Y_n^{4k}\circ \tilde{R}_n$ lying
between any two nodal meridians of $Y_n^{2k}$.  Accordingly,
$Y_n^{4k}\circ \tilde{R}_n$ takes {\em negative} values at all the points
$P_j$; and {\em positive} values at all the points
$Q_j,\; 1\leq j\leq k\cdot 4k$.
\end{lemma}

We can next adjust $\tilde{R}_n$ slightly to get a rotation $R_n$
so that the conclusion of Lemma \ref{lemma:1mod4} still holds,
and also $Y_n^{4k}\circ R_n$ takes
\begin{itemize}
\item[i)] {\em Positive} value at $S$, and {\em negative} value at $N$.
\item[ii)]  {\em Positive} values at equatorial intersection
points for which $0<\phi\leq\pi$, and {\em negative} values
at equatorial intersection points for which $\pi<\phi\leq 2\pi$
\end{itemize}

It is then clear that if we choose $\epsilon_n$ small enough,
then the nodal set of $Y_n^{2k}+\epsilon_n\cdot Y_n^{4k}\circ R_n$
will consist of
\begin{itemize}
\item[i)] one ``long'' oval along the equator.
\item[ii)] {\em at least} $k\cdot 2k$ ``small'' ovals bounding
{\em positive} nodal domains in the upper hemisphere.
\item[iii)] {\em at least} $k\cdot 2k$ ``small'' ovals bounding
{\em negative} nodal domains in the lower hemisphere.
\end{itemize}
Actually, there will be some additional small nodal domains
in upper and lower hemisphere, but we shall ignore them.

Altogether, the nodal set consists of at least $4k^2+1\sim n^2/4$
disjoint ovals, finishing the proof for $n=4k+1$.
\qed

%%%%%%%%%%%%%

{\em Case iii): $n=4k, n=4k+2$.}

Let $n=2m$.  We shall be perturbing the spherical harmonics
$Y_{2m}^m$ by (small multiples of) $Y_{2m}^{2m}$ which is
proportional to $\sin(2m\phi)$. Nodal meridians of $Y_{2m}^m$ are
the meridians
$$
\cN_1=\{\phi=0,\phi=\pi/m,\phi=2\pi/m,\ldots,\phi=\pi,\phi=\pi(m+1)/m,
\ldots\},
$$
while those of $Y_{2m}^{2m}\circ R(-\delta)$ (here $R$ denotes a
rotation around $NS$-axis; we shall chose $\delta\ll 1/2m$) have the
form
$$
\cN_2=\{\phi=\delta,\phi=\pi/(2m)+\delta,\phi=2\pi/(2m)+\delta,
\ldots\}
$$

It is easy to show that
\begin{lemma}\label{lem1:even}
For small enough $\delta$ there will be {\em two} meridians from
$\cN_2$ between any two meridians of $\cN_1$.
\end{lemma}

Accordingly, one can arrange for $Y_{2m}^{2m}\circ \tilde{R}_{2m}$
to have constant (say, positive) sign at all double intersections of
the nodal lines of $Y_{2m}^m$ different from $N$ and $S$.  By
adjusting the rotation slightly, we can arrange for
$Y_{2m}^{2m}\circ \tilde{R}_{2m}$ to be nonzero at $N$ and $S$ as
well (and to have positive sign there as well).

Accordingly, for small enough $\epsilon$, the nodal set of the
spherical harmonic $Y_{2m}^m+\epsilon Y_{2m}^{2m}\circ
\tilde{R}_{2m}$ will have $m(m+1)\sim n^2/4$ ovals as claimed,
finishing the proof of case ii) of Theorem \ref{thm:ovals}.

\begin{center}
\epsfxsize=4.0in \centerline{\epsffile{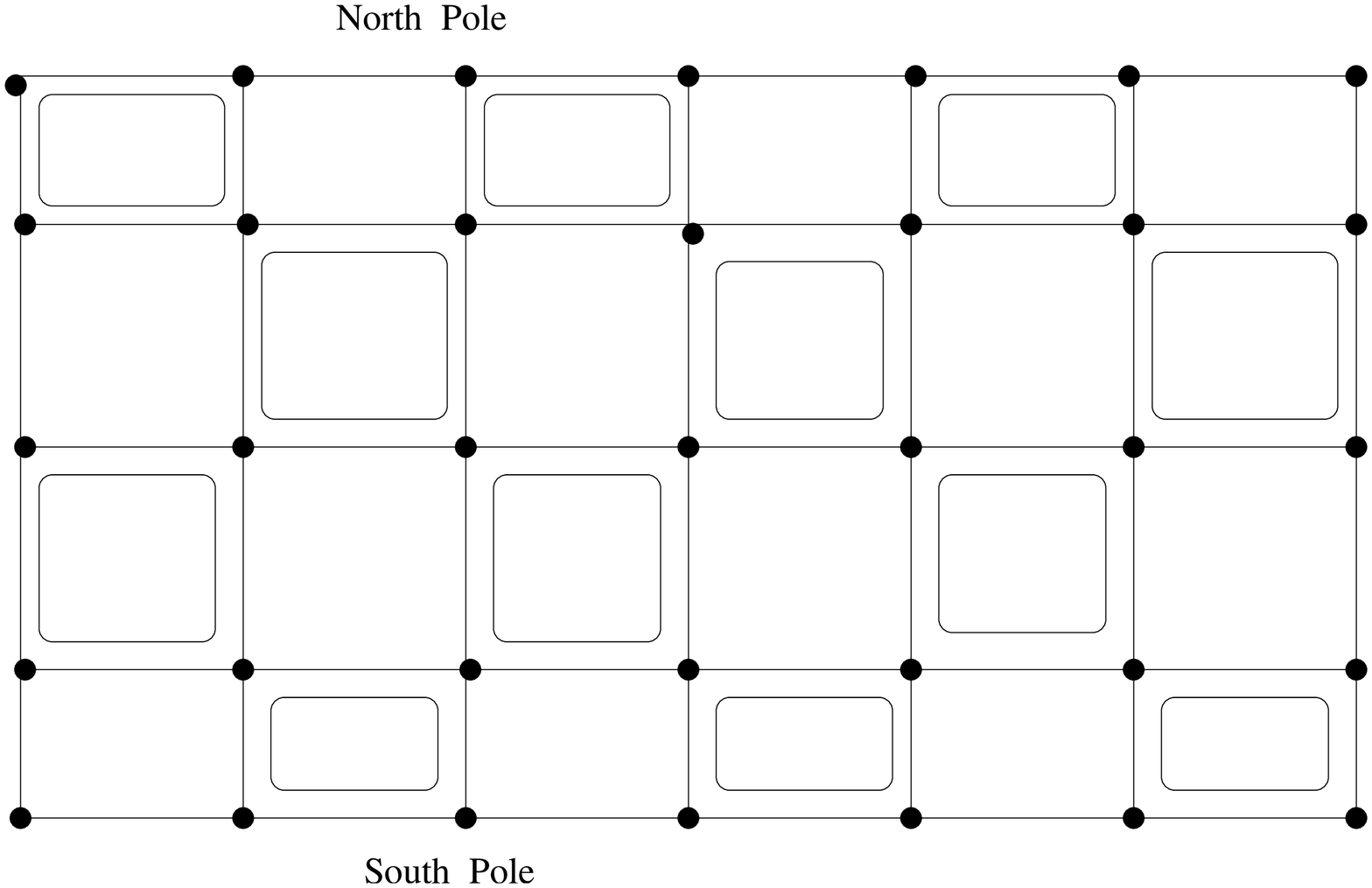}}

Figure 3. A perturbation in degree $2m=6$;\ $\bullet$\ denotes
positive sign of $Y_6^6$ at singular points of $Y_6^3$. Note that on
$S^2$ all points on the top line (north pole) and the bottom line
(south pole) are identified.
\end{center}
\qed

%%%%%%%%%%%%%

\section{Eigenfunctions in $\reals^2$ with two nodal domains}

Let $M$ be a compact surface, and let $\Delta$ be the
Laplace-Beltrami operator for a metric on $M$, and let
$\Delta u_k=\lambda_k u_k$ be a sequence of eigenfunctions.
Let $x_k\in M$ be a point of supremum of $u_k$,
and let $B_k$ be a geodesic ball
centered at $x_k$ with radius $C/\sqrt{\lambda_k}$.  Blow up
$B_k$ to the unit disk in $\reals^2$, and let
$\tilde{u}_k/\sup u_k$ be the eigenfunction after that change of
variables.  Then (after passing to a subsequence)
$\tilde{u}_k$ converges to a solution
\begin{equation}\label{eig1}
\Delta\tilde{u}=\tilde{u}, \qquad |\tilde{u}|<1.
\end{equation}

Thus the nodal structure of $\tilde{u}$  models fine local nodal
structure of $u_k$, e.g. if one is able to prove that $u$ has
infinite number of isolated critical points, then one can expect
their number to be unbounded also for the  sequence $u_k$. We guess
here that any solution of \ref{eig1} on $\reals^2$ has infinite
number of critical points. However, a bit stronger statement on the
infinite number of nodal domains of $\tilde{u}$ turns not to be
true.

Somewhat surprisingly, we prove the following
\begin{theorem}\label{2domains}
There exist solutions of \eqref{eig1} with exactly {\em two} nodal
domains.
\end{theorem}

{\em Proof of Theorem \ref{2domains}.}
The proof is inspired by a related construction in \cite{Lewy}.
Take a function $f(x,y)$ defined in polar
coordinates $x=r\cos\theta,y=r\sin\theta$ by
$$
f=J_1(r)\sin\theta.
$$
Then zeros of $f$ consist of the $x$-axis, together with
an infinite union of concentric circles
$$
r= j_k, k=1,2,\ldots
$$
where $j_k$ denotes the $k-th$ zero of $J_1$.

We shall use the following standard result (see e.g. \cite[Thm
1.6]{Segura}):
\begin{lemma}\label{lem1:bessel}
There exists $C>0$ such that
$$
\inf_k|j_{k+1}-j_k|>C.
$$
\end{lemma}

An easy consequence of Lemma \ref{lem1:bessel} is the following
\begin{lemma}\label{lem2:bessel}
Let $0<\delta_2<\delta_1< C/2$; define function $g(x,y)$ by the
formula
$$g(x,y):=f(x-\delta_1,y-\delta_2).$$
Then $g(j_k,0)=(-1)^{k-1}$ and $g(-j_k,0)=(-1)^k$.
\end{lemma}

Finally, it is easy to see that the function
$$
h(x,y):=f(x,y)+\epsilon\cdot g(x,y)
$$
will have the properties required in Theorem \ref{2domains} for
$\epsilon$ small enough.
\qed

\begin{center}
\epsfxsize=3.0in \centerline{\epsffile{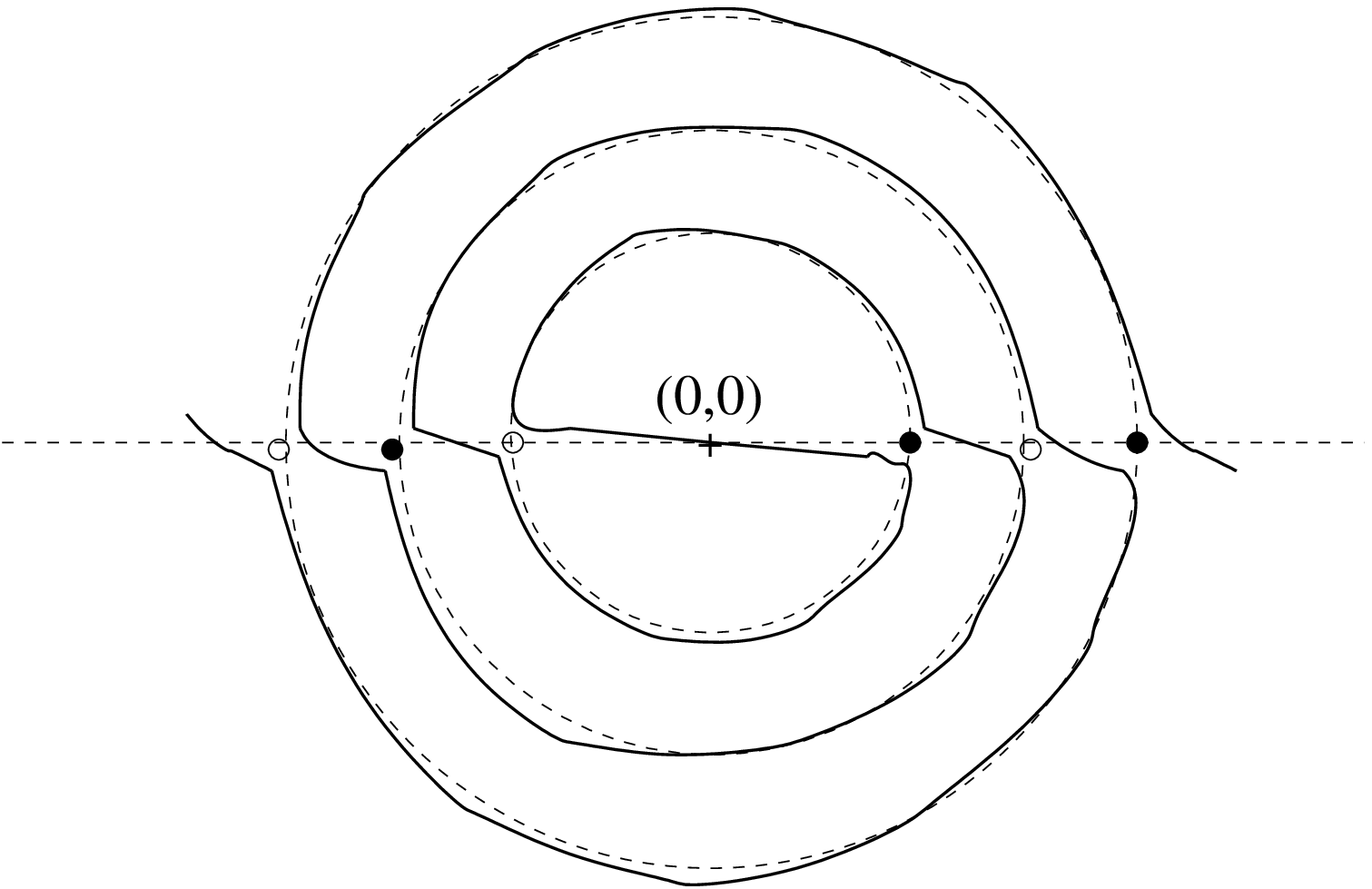}} Figure 4. A
function in $\reals^2$ with two nodal domains. Dashed line denotes
the nodal set of $f=J_1(r)\sin\theta$, solid line denotes the nodal
set of the perturbed eigenfunction $h$.\ $\bullet$\ denotes positive
sign of $g$ at singular points of $f$, \ $\circ$\ denotes negative
sign of $g$ at singular points of $f$.
\end{center}

\

\

\noindent{\bf Acknowledgements.} The authors would like to thank S.
Lando, F. Nazarov, K. Przybytkowski, M. Sodin and A. Zvonkin for
useful discussions. The authors would like to thank the anonymous
referee for many useful corrections and suggestions that helped
greatly improve the paper.

The first and second authors would like to thank the organizers of
the 2006 Spring Workshop on Dynamical Systems and Related Topics in
honor of Y. Sinai for their hospitality during that conference.  A
large part of this paper was written while the second author visited
IHES, France; their hospitality is greatly appreciated.


\begin{thebibliography}{1}

\bibitem{unsolved} V. Arnold, M. Vishik,  Y. Ilyashenko, A. Kalashnikov,
V. Kondratyev, S. Kruzhkov, E. Landis, V. Millionshchikov, O.
Oleinik, A. Filippov, M.Shubin. Some unsolved problems in the theory
of differential equations and mathematical physics. Uspekhi Mat.
Nauk 44 (1989), no. 4(268), 191--202; transl. in Russian Math.
Surveys 44 (1989), no. 4, 157--171.

\bibitem{CH} R. Courant and D. Hilbert. Methods of mathematical
physics, vol. I.  Interscience Publishers, Inc., New York, N.Y.,
1953.

\bibitem{EG} A. Eremenko and A. Gabrielov. Rational functions
with real critical points and the B. and M. Shapiro conjecture in
real algebraic geometry, Annals of Math. 155 (2002), no. 1,
105--129.

\bibitem{Gudkov} D. Gudkov.  The topology of real projective algebraic varieties.
Uspehi Mat. Nauk 29 (1974), no. 4(178), 3--79.

\bibitem{Survey} D. Jakobson, N. Nadirashvili and J. Toth.
Geometric properties of eigenfunctions, Russian Math. Surveys,
56, 6 (2001) 67--88.

\bibitem{Karp89} V. N. Karpushkin. Topology of zeros of
eigenfunctions. Funct. Anal. Appl. 23 (1989),  no. 3, 218--220 (1990).

\bibitem{Karp94} V. N. Karpushkin. The number of components of the
complement of the level surface of a harmonic polynomial in three variables.
Funct. Anal. Appl.  28  (1994),  no. 2, 116--118.

\bibitem{Karp02} V. N. Karpushkin. On the number of components of
the complement to some algebraic curves.  Russian Math. Surveys
57 (2002),  no. 6, 1228--1229.

\bibitem{LZ}  S. Lando and A. Zvonkin. Graphs on surfaces and their
applications. With an appendix by Don B. Zagier. Encyclopaedia of
Mathematical Sciences, 141. Low-Dimensional Topology, II.
Springer-Verlag, Berlin, 2004.

\bibitem{Lewy} H. Lewy, On the minimum number of domains
in which the nodal lines of spherical harmonics divide the sphere.
Comm. PDE, 2(12) (1977) 1233--1244.

\bibitem{Leydold} J. Leydold. On the number of nodal domains of spherical
harmonics. Topology 35 (1996), no. 2, 301--321.

\bibitem{generic} J. Neuheisel. Asymptotic distribution of
nodal sets on spheres. Thesis, Johns Hopkins University, Baltimore,
MD 2000,  28  (1994),  no. 2, 116--118.
\newline
http://mathnt.mat.jhu.edu/mathnew/Thesis/joshuaneuheisel.pdf

\bibitem{Pleijel} A. Pleijel. Remarks on Courant's nodal line theorem.
Comm. Pure Appl. Math 9 (1956), 543--550.

\bibitem{santos} F. Santos.   Optimal degree construction of real
algebraic plane nodal curves with prescribed topology. I. The
orientable case. Real algebraic and analytic geometry (Segovia,
1995). Rev. Mat. Univ. Complut. Madrid 10 (1997), Special Issue,
suppl., 291--310.

\bibitem{Segura} J. Segura. Bounds on differences of adjacent
zeros of Bessel functions and iterative relations between consecutive zeros.
Math. Comp. 70 (2001), no. 235, 1205--1220.

\bibitem{TS} A. Tikhonov and A. Samarskii. Equations of mathematical
physics, Moscow, 1953.

\bibitem{Viro} O. Viro, Real algebraic plane curves: constructions
with controlled topology. Leningrad Math. J., 1, 5 (1990)
1059--1134.

\end{thebibliography}
\end{document}